\documentclass{article}
\usepackage{graphicx} 

\usepackage{amssymb,amsthm,amsmath,marvosym,amsfonts,fontenc}
\usepackage{comment}
\usepackage{enumerate,marginnote}
\usepackage{graphicx}
\usepackage{hyperref}
\usepackage{cleveref}

\newtheorem{theorem}{Theorem}
\newtheorem*{theorem*}{Theorem}
\newtheorem*{observation*}{Observation}

\newtheorem{subtheorem}{Theorem}[theorem]

\newtheorem{lemma}[theorem]{Lemma}

\newtheorem{conjecture}[theorem]{Conjecture}

\newtheorem{proposition}[theorem]{Proposition}

\theoremstyle{remark}
\newtheorem{claim}[subtheorem]{Claim}

\def\epsilon{\varepsilon}

\def\phi{\varphi}

\usepackage{graphicx,color}
\title{Totally odd immersions in line graphs}

\author{Andrea Jim\'enez\thanks{A. Jiménez is partially supported by ANID/Fondecyt Regular 1220071, MATHAMSUD
  MATH210008, and ANID/PCI/REDES 190071. Email: andrea.jimenez@uv.cl} \hspace{1.4cm} Daniel A. Quiroz\thanks{D.A. Quiroz is partially supported by ANID/Fondecyt Iniciación en
  Investigación 11201251, and MATHAMSUD MATH210008. Email: daniel.quiroz@uv.cl} \\ \small{Instituto de Ingenier\'ia Matem\'atica - CIMFAV, Universidad de Valpara\'iso, Chile} \\ \\ Christopher Thraves Caro\thanks{Email: cthraves@udec.cl }\\ \small{Departamento de Ingenier\'ia Matem\'atica, Universidad de Concepci\'on, Chile}}

\date{}

\begin{document}

\maketitle

\begin{abstract}
    The immersion-analogue of Hadwiger's Conjecture states that every graph $G$ contains an immersion of $K_{\chi(G)}$. This conjecture has been recently strengthened in the following way: every graph~$G$ contains a totally odd immersion of $K_{\chi(G)}$. We prove this stronger conjecture  for line graphs of constant-multiplicity multigraphs, thus extending a result of Guyer and McDonald.
\end{abstract}

\section{Introduction}

Every graph considered in this paper is loopless and simple whenever not explicitly called a multigraph. A graph~$G$ is said to contain another graph $H$ as an \emph{immersion} if there exists an injective function $\phi\colon V(H)\rightarrow V(G)$ such that:
\begin{enumerate}[(I)]
\item For every $uv\in E(H)$, there is a path in $G$, denoted $P_{uv}$, with endpoints $\phi(u)$ and~$\phi(v)$.
\item The paths in $\{P_{uv} \mid uv\in E(H) \}$ are pairwise edge-disjoint.
\end{enumerate}
The vertices in $\phi(V(H))$ are called the \emph{terminals} of the immersion. If the terminals are not allowed to appear as interior vertices on paths $P_{uv}$, then $G$ is said to contain $H$ as a \emph{strong immersion}, and if, moreover, the paths are internally vertex-disjoint, then $G$ is said to contain a \emph{subdivision} of $H$. An immersion (or strong immersion) is said to be \emph{totally odd} if all the paths are odd.

Haj\'os conjectured that every graph $G$ contains a subdivision of a clique on $\chi(G)$ vertices. Catlin~\cite{C79} showed that this conjecture is false, his counterexamples being line graphs of constant-multiplicity multigraphs. On the other hand, in 2007, Thomassen gave an elegant two-page proof of the following, which implies that the conjecture of Haj\'os is true for line graphs of graphs.

\begin{theorem}[Thomassen~\cite{T07}]\label{thm:thomassen}
    Every (simple) graph $G$ with maximum degree $d$ and chromatic index
$d + 1$, contains two vertices $x,y$ and a collection of $d$ pairwise edge-disjoint paths between
$x$ and~$y$.
\end{theorem}

A weaker conjecture $-$ still open $-$ was proposed by Abu-Khzam and Langston \cite{AKL03}: every graph~$G$ contains an immersion of $K_{\chi(G)}$. Theorem~\ref{thm:thomassen} implies that this conjecture holds for line graphs of graphs, and Guyer and McDonald~\cite{GMcD19} used this fact to show that the conjecture holds for line graphs of constant-multiplicity multigraphs.  In particular, the examples of Catlin satisfy this weaker conjecture (see also~\cite{Q20}).

The conjecture of Abu-Khzam and Langston has recently been strengthened as follows.

\begin{conjecture}[Churchley \cite{C17}]\label{conj:church}
Every graph $G$ with $\chi(G)\ge t$ contains a totally odd immersion of $K_t$.
\end{conjecture}

Note that Conjecture~\ref{conj:church} and that of Haj\'os are incomparable, because subdivisions and totally odd immersions  are incomparable: the class of bipartite graphs contains all complete graphs as subdivisions, but excludes $K_3$ as a totally odd immersion, while the class of planar graphs contains all complete graphs as totally odd immersions\footnote{To see this take a star with, say, $t$ leaves, and replace each edge by a multiedge with multiplicity $t-1$. Subdivide each new edge once, to obtain a graph that contains $K_t$ as an immersion, and further subdivide some edges as needed to make the paths odd.}, while excluding $K_5$ as a subdivision.

As evidence for Conjecture~\ref{conj:church}, we prove the following, which extends the result of Guyer and McDonald.


\begin{theorem}\label{thm:main}
If $G$ is the line graph of a constant-multiplicity multigraph, then $G$ contains a totally odd strong immersion of $K_{\chi(G)}$.
\end{theorem}



In light of Theorem~\ref{thm:main} we make the following strengthening of Conjecture~\ref{conj:church}

\begin{conjecture}\label{conj:ours}
Every graph $G$ with $\chi(G)\ge t$ contains a totally odd strong immersion of $K_t$.
\end{conjecture}

Further evidence for Conjecture~\ref{conj:ours} can be found in \cite{BQSZ} where an analogue  of the theorem of Duchet and Meyniel \cite{DM82} is proved for totally odd strong immersions. Moreover, it is shown in \cite[Corollary 5.4]{C17} that  Conjecture~\ref{conj:ours} holds for ``almost every'' graph.

Other related results include the proof by Reed and Seymour of Hadwiger's Conjecture for line graphs (of multigraphs)~\cite{RS04},  and Steiner's recent proof of the so-called Odd Hadwiger Conjecture for line graphs of graphs~\cite{S23}. Note that minors and immersions are incomparable in general, and, moreover, Guyer and McDonald showed that they are incomparable within line graphs. 

\section{The proof}

We use the following result which is implied by Vizing's Adjacency Lemma (see e.g. \cite[p. 71]{FW95}).

\begin{lemma}[Vizing]\label{lem:vizing}
Let $H$ be an edge-critical graph with maximum degree $d$. Then every vertex of $H$ has at least two neighbors of degree $d$.
\end{lemma}

We first prove our result for line graphs of graphs.

\begin{theorem}\label{thm:strong}
Every line graph $G$ of a (simple) graph contains a totally odd strong immersion of $K_{\chi(G)}$.
\end{theorem}

\begin{proof}
 Suppose (reductio ad absurdum) that $G$ is a minimal counterexample, and let $H$ be such that $G=L(H)$. We can assume that $H$ has maximum degree $d$ while $\chi(G)=d + 1$, and, moreover, we can assume $d\ge 3$. By Theorem~\ref{thm:thomassen}, $H$ contains two vertices $x,y$ and a collection $P_1, \dots, P_d$ of pairwise edge-disjoint paths between $x$ and $y$. For $i\in [d]$, let $P_i=x,h^i_1, \dots ,h^i_{\ell_i},y$, and let $Q_i=v^i_0,\dots ,v^i_{\ell_i}$ be the path in $G$ such that $V(Q_i)=E(P_i)$. Note that we can have $v^i_0=v^i_{\ell_i}$ for at most one value of~$i$, and if so, we assume that this value is $i=d$.  

If for some $i\in [d]$ we have $\ell_i\ge 2$, then at least one of $h^i_1, h^i_2$ has degree at least 3, as otherwise one of these two vertices would not have two neighbors of degree $d$, contradicting Lemma~\ref{lem:vizing} as $H$ is edge-critical. Hence, we obtain the following.
\begin{claim}\label{claim:thirdneigh}
If $\ell_i\ge 2$, then there exists $a^i\in V(G)$ which corresponds to an edge $e^i \notin E(P_i)$ incident to $h^i_1$ or to $h^i_2$.
\end{claim}




To bring about the desired contradiction, we will construct a totally odd immersion of $K_{d+1}$. We let $G[v^1_0,v^2_0\dots , v^d_0]$ together with some $v^*\notin \{v^1_0,\dots ,v^d_0\}$ be part of that immersion, where $v^*= v^1_{\ell_1}$ unless otherwise stated. 
We are only left to find the odd paths $Q_1', \dots ,Q_d'$ joining $v^*$ to every $v^i_0$. 

Let $j$ be the number of values of $i\in [d]$ such that $\ell_i$ is odd. We first consider the case in which $j \ge 4$. We can assume $\ell_1, \dots ,\ell_{j}$  are odd, and take $Q'_1=Q_1$,
$Q_i'=Q_i+v^i_{\ell_i}v^{i+1}_{\ell_{i+1}}+v^{i+1}_{\ell_{i+1}}v^1_{\ell_1}$ for $i\in \{2,\dots , j-1\}$, 
$Q_{j}'=Q_{j}+v^{j}_{\ell_{j}}v^2_{\ell_2}+v^2_{\ell_2}v^1_{\ell_1}$, and 
$Q_{i}'=Q_{i}+v^{i}_{\ell_{i}}v^1_{\ell_1}$ for $i> j$.




Assume now that $j=3$, and that $\ell_1, \ell_2 ,\ell_{3}$ are odd. If $\ell_4 \geq 2$, choose $v^*=v^4_{\ell_4}$ and take 
$Q_i'=Q_i+v^i_{\ell_i}v^{i+1}_{\ell_{i+1}}+v^{i+1}_{\ell_{i+1}}v^4_{\ell_4}$ for $i\in \{1,2\}$, $Q_3'=Q_3+v^3_{\ell_3}v^{1}_{\ell_{1}}+v^{1}_{\ell_{1}}v^4_{\ell_4}$, $Q_4'=v^4_0,a^4,v^4_1,v^4_2, \dots ,v^4_{\ell_4}$ if $e^4$ is incident to $h^4_1$ or  $Q_4'=v^4_0,v^4_1,a^4,v^4_2, \dots ,v^4_{\ell_4}$ otherwise (using Claim \ref{claim:thirdneigh}), and $Q_i'=Q_i+v^i_{\ell_i}v^{4}_{\ell_{4}}$ for $i\ge 5$. If $\ell_4 = 0$, we can assume $d=4$. We can further assume that $h_1^3$ is a vertex of degree at least 3 in $H$, due to Lemma \ref{lem:vizing}. So, there exists a vertex $a$ in $G$ which corresponds to an edge not in $P_3$ incident to $h^1_3$. We take $Q_1'=Q_1$, $Q_2'=Q_2 +v^2_{\ell_2}v^3_{\ell_3}+v^3_{\ell_3}v^1_{\ell_1}$, 
$Q_3'=v^3_0,a,v^3_1, \dots ,v^3_{\ell_3}  
+v^3_{\ell_3}v^4_{\ell_4}+v^4_{\ell_4}v^2_{\ell_2}+v^2_{\ell_2}v^1_{\ell_1}$,  and
$Q_4'=Q_4+v^4_{\ell_4}v^1_{\ell_1}$. Thus, for $j=3$, we are only left with the case $d=3$. We can assume that at least two of $h^1_1,h^2_1, h^3_1$ have degree 3, say  $h^2_1$ and $h^3_1$, as otherwise $x$ would contradict Lemma~\ref{lem:vizing}. So there exist vertices $b^2, b^3$ in $G$ which correspond to edges incident to $h^2_1, h^3_1$ and not in $E(P_1 \cup P_2 \cup P_3)$. We take $Q_1'=Q_1$, $Q_2'=v^2_0,b^2,v^2_1,v^2_2\dots ,v^2_{\ell_2}+v^2_{\ell_2}v^1_{\ell_1}$,
and $Q_3'=v^3_0,b^3,v^3_1,v^3_2\dots ,v^3_{\ell_3}+v^3_{\ell_3}v^1_{\ell_1}$. Note that $b^2=b^3$ is possible, but even in that case the paths remain edge-disjoint.

We consider the remaining cases. If $j=0$, then we can assume $\ell_1\ge 2$, and
take $Q_i'=Q_i+v^i_{\ell_i}v^1_{\ell_1}$ for $i\in \{2,\dots , d\}$, and $Q_1'=v^1_0,a^1,v^1_1,v^1_2, \dots ,v^1_{\ell_1}$ if $e^1$ is incident to $h^1_1$ or  $Q_1'=v^1_0,v^1_1,a^1,v^1_2, \dots ,v^1_{\ell_1}$ otherwise (using Claim \ref{claim:thirdneigh}). If $j=1$, we can assume $\ell_1$ is odd, and take $Q_1'=Q_1$ and $Q_i'=Q_i+v^i_{\ell_i}v^1_{\ell_1}$ for $i\in \{2,3, \dots ,d\}$. So we are only left with the case $j=2$.
If $d\geq 4$, we take $Q_1'=Q_1$, $Q_2'=Q_2 +v^2_{\ell_2}v^3_{\ell_3}+v^3_{\ell_3}v^1_{\ell_1}$, 
$Q_3'=Q_3+v^3_{\ell_3}v^4_{\ell_4}+v^4_{\ell_4}v^2_{\ell_2}+v^2_{\ell_2}v^1_{\ell_1}$,  and
$Q_i'=Q_i+v^i_{\ell_i}v^1_{\ell_1}$ for $i\in \{4,\dots , d\}$. 
If $d=3$, we can assume that at least one of $h^1_1,h^2_1$ has degree 3, say $h^2_1$, as otherwise $x$ would contradict Lemma~\ref{lem:vizing}. So there exists a vertex $b$ in $G$ which corresponds to the edge incident to $h^2_1$ and not in $P_2$. We take $Q_1'=Q_1$, $Q_2'=v^2_0,b,v^2_1,v^2_2\dots ,v^2_{\ell_2}+v^2_{\ell_2}v^1_{\ell_1}$,
and $Q_i'=Q_i+v^i_{\ell_i}v^1_{\ell_1}$ for $i\in \{3,\dots , d\}$. 










In all these cases, we obtain a totally odd strong immersion of $K_{d+1}$ in $G$, a contradiction.
\end{proof}

We need a lemma to extend Theorem~\ref{thm:strong} to line graphs of constant-multiplicity multigraphs, and, for this, notation. Let $G$ be a graph and $B_m(G)$ be the \emph{$m$-blow-up} of $G,$ that is, we replace each vertex $x$ of $G$ by an independent set $x^1,\dots, x^m$, and we let $x^1,\dots, x^m$ and $y^1,\dots,y^m$  induce a complete bipartite graph whenever $xy\in E(G)$. The next lemma can be deduced from the proof of \cite[Theorem 1.2]{GMcD19}.

\begin{lemma}[Guyer and McDonald~\cite{GMcD19}]\label{lemma:jessica}
Let $P=u_0,u_1, \dots , u_\ell$ be a path. The graph $B_m(P)$ contains, for each pair $1\le i,j\le m$, a path of the form $u_0^{k_0},u_1^{k_1}, \dots u^{k_\ell}_\ell$, with $k_0,\dots ,k_\ell\in [m]$, where  these paths are pairwise edge-disjoint.
\end{lemma}

For a multigraph $H$, and an integer $m\ge 2$, we let $mH$ be the multigraph obtained from $H$ by replacing each edge with a parallel edge of multiplicity $m$. Theorem~\ref{thm:main} follows from Theorem~\ref{thm:strong} and the following result, because if $\chi(L(H))=t$, then it is easy to see that $\chi(L(mH))\le mt$. 

\begin{theorem}
    Let $H$ be a multigraph such that $L(H)$ contains a totally odd strong immersion $K_t$, and $m\ge 2$ be an integer. Then $L(mH)$ contains a totally odd strong immersion of $K_{mt}$.
\end{theorem}

\begin{proof}
Let $G=L(H)$, and let $S\subseteq V(G)$ be the set of terminals of a  totally odd immersion of $K_t$ in~$G$. Further, let $\mathcal{P}$ be the set of paths of that immersion, and $S_m$ be the set that contains the $m$ copies of~$v$ in $L(mH)$ for every $v\in V(G)$. For every $x,y$-path $P$ of $\mathcal{P}$, Lemma~\ref{lemma:jessica} tells us that there is a set of pairwise edge-disjoint paths that joins each copy of $x$ to each copy of $y$. Since the path $P$ is odd, the lemma guarantees that all these new paths are odd. Moreover, since the paths in $\mathcal{P}$ are pairwise edge-disjoint, if $P_1,P_2\in \mathcal{P}$, then any path obtained from $P_1$ using the lemma will be edge-disjoint from any path obtained from $P_2$ using the lemma. We conclude that $L(mH)$ contains an immersion of $K_{mt}$ with $S_m$ as its set of terminals.
\end{proof}





\end{document}